\documentclass{jams-l}

\newcommand{\Q}{\Bbb{Q}}
\newcommand{\ctd}{t_\Psi}
\renewcommand{\span}{\operatorname{span}}
\renewcommand{\max}{\operatorname{max}}
\newcommand{\dottedrightarrow}{\rightarrow}
\renewcommand{\div}{\operatorname{div}}
\renewcommand{\bullet}{*}

\newcommand{\lvol}{\operatorname{vol}}
\newcommand{\Td}{\operatorname{Td}}
\newcommand{\Spec}{\operatorname{Spec}}
\newcommand{\sigone}{\Sigma ^{(1)}}
\newcommand{\Cal}{\mathcal}

\newcommand{\mult}{\text{mult}\,}

\newcommand{\lp}{\langle}
\newcommand{\rp}{\rangle}

\newcommand{\Hom}{\operatorname{Hom}}
\newcommand{\tupler}{\rho_1^{(a_1)},\dots,\rho_r^{(a_r)}}

\theoremstyle{plain}
\newtheorem{theorem}{Theorem}
\newtheorem{corollary}{Corollary}
\newtheorem{proposition}{Proposition}
\newtheorem{lemma}{Lemma}
\newtheorem*{kludge}{Weak Oda Theorem \cite{W}}

\theoremstyle{definition}
\newtheorem*{definition}{Definition}
\newtheorem{example}{Example}

\begin{document}
\title {Cycles representing the Todd class of a toric variety}

\begin{abstract}
In this paper, we describe a way to construct cycles which represent the
Todd class of a toric variety.  
Given a lattice with an inner product,
we assign a rational number $\mu(\sigma)$ 
to each
rational polyhedral cone $\sigma$ in the lattice, 
such that for any toric variety
$X$ with fan $\Sigma$ in the lattice, we have:

$$\operatorname{Td}(X)=\sum_{\sigma \in \Sigma} \mu(\sigma) [V(\sigma)].$$
This constitutes an improved answer to an old question of Danilov 
\cite{Dan}.  

In a similar way, beginning with the choice of a complete flag in the lattice,
we obtain the cycle Todd classes constructed by Morelli in 
\cite{Mor}.

Our construction is based on an intersection
product on cycles of a simplicial toric variety developed by the second-named
author in \cite{Tho}.  Important properties of the construction are established
by showing a connection to the canonical representation of the Todd class
of a simplicial toric variety as a product of torus-invariant divisors
developed by the first-named author in \cite{Pom2}.
\end{abstract}

\author
{James Pommersheim}
\address
{ Department of Mathematics,
Pomona College,
Claremont, California, 
92037 USA }
\email{jpommersheim@pomona.edu}

\author{Hugh Thomas}
\address{Fields Institute,
222 College Street,
Toronto ON,
M5T 3J1 Canada}
\email{hthomas@fields.utoronto.ca}



\maketitle

\section{Introduction}

Shortly after toric varieties were introduced in the 1970's, Danilov
\cite{Dan} 
asked whether it
is possible, given a lattice $N$, to assign a rational number $\mu(\sigma)$
 to every cone $\sigma$ in $N$ so that for any fan $\Sigma$ in
the lattice, the Todd class $\Td(X)$ of the corresponding toric variety $X$
would be given by 
\begin{equation}
\operatorname{Td}(X) = \sum_{\sigma\in \Sigma} \mu(\sigma) [V(\sigma)],
\end{equation}
where $[V(\sigma)]$ is the class in the Chow group of $X$ corresponding
to the orbit closure $V(\sigma)$.
Note that the above formula is to hold simultaneously for all fans in $N$.  
(The much weaker statement
that the Todd class of
any fixed toric variety can be written as some linear combination of the
$[V(\sigma)]$ follows from the fact that the $[V(\sigma)]$ generate the Chow
group of $X$.)

Morelli answered Danilov's question by showing (non-constructively) 
the existence of
such an assignment \cite[Theorem 3]{Mor}. He also gave an explicit
canonical construction of such an assignment, where the $\mu(\sigma)$
are rational functions on a certain Grassmannian \cite[Theorem 11]{Mor}.
In this paper, we construct  assignments satisfying Danilov's conditions,
in which the coefficients $\mu(\sigma)$ are rational
 numbers (Corollary 1).

The present work also
serves to unify two previous formulas for the Todd class of a toric
variety: Morelli's canonical expression for the Todd class in terms
of rational functions on Grassmannians   \cite[Theorem 11]{Mor},
and the first-named author's expression for the Todd class of
a simplicial toric variety as a polynomial in the torus-invariant
divisors \cite{Pom2}.  In particular, we show that if this polynomial
expression is multiplied out using the ring structure
developed in \cite{Tho}, we obtain 
Morelli's canonical formula (Corollary 2).

One motivation for Danilov's question is the connection between Todd classes 
of toric varieties and counting lattice points in polytopes.  
We will now briefly sketch the application of our results to this problem.  
If $P$ is a lattice polytope in $\Bbb{Z}^n$, and 
$F$ is a face of $P$,
there
is an associated outer normal cone, which we denote by $C(F,P)$. 
 Also associated to $F$ is 
a lattice volume, $\lvol(F)$: this is the volume of $F$ normalized
with respect to the lattice in the affine span of $F$.  
In this paper, we construct a map $\mu$ from rational cones in 
$\Bbb{Z}^n$ to $\Q$ with the following property:  for any lattice polytope $P$ 
in $\Bbb{Z}^n$, 
the number of lattice points in $P$ is given by the formula

\begin{equation} \#(P)=\sum_{F} \mu(C(F,P)) \lvol(F),
\end{equation}
where the sum is taken over all faces $F$ of $P$.

The existence of maps $\mu$ satisfying (2) was proved by
McMullen [McM], but his proof is not constructive.  The difficulty is
that (2) does not determine $\mu$ uniquely, so one must make some 
choices.  
The paper [BP] shows a way of unwinding  McMullen's proof to yield
what is in principle a constructive definition of a function $\mu$ satisfying
(2), based on a 
initial choice of an ordered basis for the lattice containing $P$.  
However, computing 
their $\mu$ for any given cone is a rather complicated procedure.

\medskip

We now describe our construction of a cycle Todd class 
in more detail. Let us suppose that we have a simplicial toric variety 
$X$ defined by a fan
$\Sigma$ in a lattice $N$.  
We recall from \cite{Tho} a construction which shows 
that, having selected a {\it complement map} 
for 
$\Sigma$, we obtain a ring structure on $Z_*(X)_\Q$, the torus-invariant 
cycles of 
$X$
with rational coefficients.  
(We will define the notion of a complement map later; for now, let it suffice
to say that choosing an inner product on $N\otimes \Bbb{Q}$ is one way to
determine
a complement map for any fan $\Sigma$ in $N$.)  

$A_*(X)_\Q$, the Chow ring of cycles of 
$X$ modulo rational equivalence,  is naturally a quotient 
of $Z_*(X)_\Q$;
the ring structure on $Z_*(X)_\Q$ is defined so that this quotient map
is a ring homomorphism. 

There is a well-known expression for the Todd class of a smooth toric 
variety $X$ 
as a 
polynomial in the torus-invariant divisors, evaluated in $A_*(X)_\Q$.
Once we have chosen a complement map for the lattice, we may evaluate
 this polynomial in 
$Z_*(X)_\Q$ instead, and obtain a cycle that
 represents the Todd class of $X$. The 
formula 
for the dimension $k$ part of this cycle Todd class 
consists of a sum over cones 
of 
$\Sigma$ 
of codimension $k$ of a contribution which depends only on the cone itself, and
not on data from other parts of the fan.  

If our toric variety $X$ is not smooth, the standard way to compute the Todd 
class 
is to take a resolution of singularities $\pi:X' 
\rightarrow X$, 
and then push forward the Todd class of 
$X'$.  We wish to mimic this process to construct a cycle Todd class for $X$.  
There is no problem taking a toric resolution of singularities, and the 
construction described above gives us a cycle representing the Todd class of 
$X'$, 
which we can push forward to $X$.  We denote this cycle by $\ctd^{X'}(X)$,
where $\Psi$ denotes the complement map.

The question that naturally arises is whether  $\ctd^{X'}(X)$  depends
 on $X'$. Theorem 1 of this paper asserts
 that for a certain class of complement maps (called
{\it rigid} complement maps, which include those induced from a choice of
inner product),   $\ctd^{X'}(X)$ is independent of 
$X'$.  It follows that for any such $\Psi$, we have
constructed an assignment of rational numbers to cones that satisfies
Danilov's condition  (Corollary 1). 

Our proof of Theorem 1 relies on 
the formula of \cite{Pom2} which expresses the Todd
class of a complete, 
simplicial toric variety in a canonical way as a product of 
torus-invariant divisors
with coefficients determined by the local combinatorics of the fan.
Specifically, a function $f$ is constructed from tuples $\tupler$
of rays in a lattice $N$ with multiplicities $a_i>0$ to the rational
numbers which satisfies
$$
\operatorname{Td}^k(X) = \sum f(\tupler)[V(\rho_1)]^{a_1}\cdots [V(\rho_r)]^{a_r}
$$
where the sum is taken of all tuples of rays in $\sigone$ with $\sum a_i
=k$.
In \cite{Pom2}, $f$ is seen to satisfy reciprocity relations such as 
those satisfied by higher dimensional Dedekind sums, and indeed the 
two-dimensional $f(\rho_1^{(1)},\rho_2^{(1)})$ 
is seen in \cite{Pom1}
to be given by the classical Dedekind sum.  
If $X$ is complete and simplicial, then given a complement map
for $X$, we obtain a cycle Todd class for $X$ by evaluating the formula of 
\cite{Pom2} in $Z_*(X)_\Q$.  We denote this cycle Todd class by $\ctd^s(X)$.  

The connection between these two approaches to cycle Todd classes is as
follows:
we show that if
$X$ is complete and simplicial and $\Psi$ is rigid, then  
$\ctd^s(X)$ and 
$\ctd^{X'}(X)$ coincide. It follows that for any $X$, $\ctd^{X'}(X)$ does not depend on 
$X'$.

One drawback to the cycle Todd classes described above is that they
depend on a choice of complement map.  We end the paper by discussing how
to construct canonical cycle Todd classes, at the expense of expanding
the coefficient field, as is done in \cite{Mor}.

\section{The ring structure on $Z_*(X)_\Q$}

We begin by reviewing  the results of \cite{Tho}.     
Let $X$ be the toric variety determined by a fan $\Sigma$ in the lattice
$N$.  
Let $Z_i(X)$ denote the free abelian group on the torus-invariant
closed irreducible subvarieties of $X$ of dimension $i$.  
Let $Z_*(X)_\Q= \bigoplus_{i=0}^n Z_i(X)\otimes \Q$.  
(A subscript $\Q$ will always indicate tensoring with $\Q$.)

For cones $\sigma, \tau$  in $N$, let us write $\tau \prec \sigma$ if
$\tau$ is a face of $\sigma$.  When we refer to cones, we always mean
strictly convex rational polyhedral cones.  

\begin{definition}
A {\it complement map} is a map $\Psi$ from 
a set of cones
in $N$ to subspaces of $N$, satisfying the following two properties:
\begin{enumerate}
\item[(i)] $\Psi(\sigma)$ is complementary to the span of $\sigma$ in $N_\Q$.

\item[(ii)] If $\tau \prec \sigma$, and $\sigma$ is in the domain of $\Psi$, then
$\tau$ is also in the domain of $\Psi$ and $\Psi(\sigma)\subset \Psi(\tau)$.
\end{enumerate}

A complement map $\Psi$ is {\it rigid} if whenever 
$\sigma$ and $\tau$ are
two cones with the same linear span such that $\Psi(\sigma)$ is defined,
then $\Psi(\tau)$ is defined and equals $\Psi(\sigma)$. \end{definition}

There are two easy ways to construct complement maps.  One way is to put an
inner product on $N_\Q$ and let $\Psi(\sigma)$ be the orthogonal
complement of the span of $\sigma$.  Another way is to choose a flag 
$F_0\subset F_1 \subset \dots \subset F_n$ in 
$N_\Q$ and let $\Psi(\sigma)=F_k$ for all $\sigma$ of codimension $k$.  
This complement map can only be defined on those cones which intersect the
fan $F_\bullet$ generically.  
Observe that both these types of 
complement map are rigid.

The chief result of \cite{Tho} is the following:

\begin{proposition} [\cite{Tho}] If $\Sigma$ is a simplicial fan, and
$\Psi$ is a complement map 
defined on all the cones in $\Sigma$, then $\Psi$ induces 
a ring structure on $Z_*(X)_\Q$.  \end{proposition}

We wish to describe this ring structure briefly.  To do this, 
we shall have to recall
some basic notions about toric varieties.  For further details, see 
\cite{Ful}.  

Let $M$ be the dual lattice to $N$.  For $\sigma$ a rational polyhedral 
cone in $N_\Q$, let 
$\sigma\check{}$ denote the dual cone in $M_\Q$: 
$\sigma\check{}=\{x\in M_\Q|\langle x,y\rangle \geq
0$ for all $y\in \sigma\}$. Now $S_\sigma=M\cap \sigma\check{}$ is a semigroup, 
so we can form
the semigroup ring $\Bbb{C}[S_\sigma]$, writing $\chi^u$ for the basis 
element of $\Bbb{C}[S_\sigma]$ corresponding to $u\in S_\sigma$.  (The 
choice of base field is not important for our purposes; we use $\Bbb{C}$ for
sake of definiteness.)  Let $U_\sigma=\Spec(\Bbb{C}[S_\sigma])$.    

Recall that $X$ has a covering by $U_\sigma$ for 
$\sigma\in\Sigma^{\max}$, the maximal cones in $\Sigma$.
For the purposes of this exposition, we assume that all the maximal cones in 
$\Sigma$ are full-dimensional, though this is not in fact necessary.    
If $D$ is an equivariant Cartier divisor, it has, for
each $\sigma \in \Sigma^{\max}$,
a local equation on $U_\sigma$, which can be taken to be of the form
$\chi^{m_\sigma}$ for some $m_\sigma \in M$.  
This collection
of lattice points $m_\sigma$ will be referred to (somewhat abusively) as
the local equations for $D$. 

The closed subvariety $V(\tau)$ is
defined by a fan in $N/\span(\tau)$, whose cones are 
$(\sigma +\span(\tau))/\span(\tau)$ for $\sigma \succ \tau$. 
 We can write $M_\Q= \tau^\perp \oplus \Psi(\tau)^\perp$, which induces 
a  map $\pi_\tau:M_\Q \rightarrow \tau^\perp$.  
Given a Cartier divisor $D$, with local equations $m_\sigma$ for $\sigma \in 
\Sigma^{\max}$, define a Cartier divisor $D_\tau$ on $V(\tau)$ by taking its 
local equation on $(\sigma + \span(\tau))/\span(\tau)$ to be 
$\pi_\tau(m_\sigma)$ for $\sigma \in \Sigma^{\max}$.  

Now define $D\cdot V(\tau)=[D_\tau]$.  
Since any $[V(\sigma)]$ can be written as a product of $\Q$-Cartier 
divisors, this multiplication 
extends to define a multiplication on $Z_*(X)_\Q$.  (Of course, 
one has to check that this multiplication is well-defined.)

In the case where the complement map $\Psi$ is induced from an inner 
product or a complete flag, we will give a presentation for 
$Z_*(X)_\Q$. For the case where $\Psi$ comes from an inner product, 
such a presentation was given in \cite{Tho}, which we restate as Proposition 
2 below.  For the case where $\Psi$ comes from a complete flag, we give
a presentation for $Z_*(X)_\Q$ in Proposition 3 below, which is new.

 Let $\Sigma$ be a simplicial fan,
 and let $\rho_1,\dots,\rho_r$ be the rays of $\Sigma$, and
let $D_i$ be the $\Bbb{Q}$-Cartier divisor corresponding to $V(\rho_i)$. 
Let $I_\Sigma \subset \Q[D_1,\dots,D_r]$
denote the Stanley-Reisner ideal of the fan $\Sigma$, that is, the ideal
generated by products of $D_i$ whose corresponding rays are not all contained
in any cone of $\Sigma$. Let $v_i$ denote the first lattice point along
$\rho_i$.  

If the complement map is induced from an inner product, we have the
following presentation for $Z_*(X)_\Q$:

\begin{proposition} [\cite{Tho}]  For $X$ a simplicial toric variety, 
if $\Psi$ is induced from an inner product, 
then:

$$Z_*(X)_\Q \cong \Q[D_1,\dots,D_r]/I_\Sigma + \langle D_i
\sum_{1\leq j \leq r}
D_j \langle 
v_i,v_j \rangle \mid 1 \leq i \leq r \rangle$$
\end{proposition}

We now provide a similar presentation for $Z_*(X)_\Q$ in the case where
$\Psi$ is induced from a complete fan in $N_\Q$.  
   
Let us write $\div(\chi^u)$ for the principal Cartier divisor
whose local equation is $u$ on every $U_\sigma$.  Fix a complete flag
$F_\bullet$ in $N_\Q$.  
Let $J$ denote the ideal of $\Q[D_1,\dots,D_r]$ generated by terms of 
the form $\div(\chi^u)W$ where $W$ is homogeneous of degree $i$,
and $u\in F_{n-i}^\perp \cap M$.

\begin{proposition} For $X$ a simplicial toric variety, 
if $\Psi$ is induced from a complete flag 
$F_\bullet$ in
$N_\Q$,
then:

$$Z_*(X)_\Q \cong \Q[D_1,\dots,D_r]/I_\Sigma + J$$ \end{proposition}

\begin{proof} It is shown in \cite{Tho} that (no matter how $\Psi$
is chosen)
there is a surjection from
$\Q[D_1,\dots,D_r]$ to $Z_*(X)_\Q$, whose kernel contains $I_\Sigma$.   

Next, we show that if $\sigma$ is a cone of $\Sigma$ of codimension $i$, and 
$u \in F_{n-i}^\perp \cap M$, then $\div(\chi^u)\cdot [V(\sigma)]=0$.
This follows by the definition of the multiplication given above:
$\div(\chi^u)_\sigma$ is the divisor with all local equations zero.  
This shows that $Z_*(X)_\Q$ is a quotient of $\Q[D_1,\dots,D_r]/I_\Sigma+J$.
  
To complete the proof, we will show that any term in $\Q[D_1,\dots,D_r]$ is 
equivalent mod $I_\Sigma+J$ to a sum of squarefree monomials each of which
is the product of divisors corresponding to rays which span a cone of 
$\Sigma$.  This will show that, as abelian groups, $\Q[D_1,\dots,D_r]/
I_\Sigma+J \cong Z_*(X)_\Q$, which will prove the proposition.  

Suppose that we have a monomial $Q$ in $D_1,\dots,D_r$.
For the sake of simplicity, suppose that $D_1$ occurs with multiplicity
greater than 1, and the $D_i$ which occur in $Q$ are $D_1,\dots,D_s$.
So let $Q=D_1^{a_1}\dots D_s^{a_s}$  
If $\rho_1,\dots,\rho_s$ are not all contained in some cone of $\Sigma$, then
$Q$ is already contained in $I_\Sigma$.  So assume that they 
generate some cone $\sigma$.  Then,
by the genericity of $F_\bullet$, 
$F_{n-s} \cap \span(\rho_1,\dots,\rho_s)
=\{0\}$.  Thus, we can find 
$u\in F_{n-s}^\perp\cap M$, such that $\langle u,v_1\rangle \ne 0$, 
$\langle u,v_i \rangle =0$ for $2 \leq i \leq s$.  Then $\div(\chi^u)$ 
expressed as a sum of $D_i$ has $D_1$ appearing with a non-zero coefficient,
and possibly some $D_i$ with $i>s$.  Now $\div(\chi^u)\cdot[V(\sigma)] \in J$.
Thus, $\div(\chi^u)D_1^{a_1-1}D_2^{a_2}\dots D_s^{a_s}\in J$.
When we expand $\div(\chi^u)D_1^{a_1-1}D_2^{a_2}\dots D_s^{a_s}$
we obtain a non-zero multiple of $Q$  plus terms which have strictly more
distinct factors.  
Thus, mod $I_\Sigma +J$, $Q$ is
equivalent to a sum of monomials which have strictly more distinct factors.
Repeating this argument, we arrive eventually at a squarefree expression,
which completes the proof that $Z_*(X)_\Q \cong \Q[D_1,\dots,D_r]/
I_\Sigma+J$. 
\end{proof}

We now recall a couple of basic notions of intersection theory as they apply
to 
toric 
varieties.
Let $X$ and $X'$ be toric varieties defined by fans $\Sigma$ and $\Sigma'$ 
in lattices $N$ and $N'$. Let $f:X'\rightarrow X$ be a torus-equivariant
map, which implies that it is induced by a lattice map, also denoted $f$, from
$N'$ to $N$, 
such that 
for any $\sigma' \in \Sigma'$, $f(\sigma')$ is contained
in some cone of $\Sigma$.  Dualizing, $f$ also induces a map
$f^*:M\rightarrow M'$.  
  
If $D$ is a Cartier divisor on $X$ with local equations $m_\sigma$, it 
can be pulled back in the usual way to $f^*(D)$ on $X'$.  The local 
equation of $f^*(D)$ on a cone $\sigma'\in\Sigma'{}^{\text{max}}$ is 
$f^*(m_\sigma)$,
for $\sigma$ any cone of $\Sigma$ that contains $f(\sigma')$.  This notion
of pull-back clearly extends to $\Q$-Cartier divisors.  

We shall also require 
the usual push-forward map of cycles: namely, for $f:X'\rightarrow 
X$ a morphism of algebraic varieties, $W'$ a closed subvariety of $X'$, and
$W$ the Zariski closure of $f(W')$, 
$f_*([W'])=[k(W'):f^*(k(W))][W]$, if  
$W$ is of the same dimension as $W'$, and 0 otherwise.  The push-forward map
also has a natural description in terms of fans for torus-equivariant maps of
toric varieties.  This description is particularly simple in the case where 
$f:X'\rightarrow
X$ is also birational, which is the only case we shall need.
The birationality of $f$ means that we can identify $N$ and $N'$.  Let 
$\sigma'\in\Sigma'$, and let $\sigma$ be the smallest cone of $\Sigma$
containing $\sigma'$.  Then
$f_*([V(\sigma')])=[V(\sigma)]$ if $\sigma$ and $\sigma'$ have the same 
dimension, and $f_*([V(\sigma')])=0$ otherwise.  

We can now state a lemma which we shall need 
from \cite{Tho}, a cycle level version of the projection formula from 
intersection theory.  

\begin{lemma} [\cite{Tho}] 
Let $Y$ and $Z$ be simplicial toric varieties of the same dimension,  
$f:Y \rightarrow Z$ a proper torus-equivariant 
map, and $\Psi$  a rigid complement map for
both $Y$ and $Z$.  Let $D$ be a $\Q$-Cartier divisor on 
$Z$, $y\in Z_*(Y)_\Q$.  Then
$$ f_*(f^*(D)\cdot y)=D \cdot f_*(y).$$ \end{lemma}

\section{A first approach to the Todd class of a toric variety}

If $X$ is smooth, there is a formula for the Todd class of $X$, namely:
\begin{equation}
\prod_{1\leq i \leq r} \frac{D_i}{1-\exp(-D_i)}\end{equation}
where the product is taken in $A_*(X)_\Q$.  
However, this product can just as well be interpreted in our 
ring $Z_*(X)_\Q$, and the result is a cycle which represents the Todd class.  
We denote this cycle by $\ctd(X)$.  

\begin{example} Let $\{e_i\}$ be the usual basis for $\Bbb{Z}^{n+1}$.  
Define a sublattice of $\Bbb{Z}^{n+1}$ by 
$N = \{ (a_1,\dots,a_{n+1})\mid a_i \in 
\Bbb{Z},
\sum a_i=0\}$.  Let $v_i=e_i-e_{i+1}$, and $v_{n+1}=e_{n+1}-e_1$.  
Let $\Sigma$ be
the fan whose cones consist of positive linear combinations of proper subsets
of the $v_i$.  The toric variety associated to $\Sigma$ is $n$-dimensional
projective space. 

Now, put an inner product on $N$ by restriction of the usual inner product
on $\Bbb{Z}^{n+1}$, i.e. set $\langle e_i,e_j \rangle=\delta_{ij}$.  Let 
$\Psi$ be the complement map induced by this inner product.  Then the 
coefficient of a cone $\sigma \in \Sigma$ in the cycle Todd class
induced by $\Psi$ is the fraction of the linear span of $\sigma$ which is
contained in $\sigma$.  To prove this, one first reduces
to the case where $\sigma$ is 
$n$-dimensional.
By symmetry, the coefficients of all $n+1$ of the $n$-dimensional cones of $\Sigma$
must be the same.  
The sum of all the coefficients must be 1, since the (Chow group) Todd 
class has $\operatorname{Td}^n(\Bbb{P}^n)=[pt]$.  Thus the coefficients 
of each of the $n$-dimensional cones of $\Sigma$ are 
$\frac 1{n+1}$, as desired.  For further details of this calculation, see
\cite{Tho}.  
\end{example}

Smoothness for toric varieties has a simple interpretation in terms of
the fan.  We say that a cone is non-singular if the first lattice points
along its extreme rays form a basis for the lattice contained in the
linear span of the cone.  Then a toric variety $X$ is smooth iff 
all the cones of its fan are non-singular.  

Suppose $X$ is not smooth.  Then the cones of $\Sigma$ can be subdivided, 
yielding a fan $\Sigma'$, in such a way that all the cones of $\Sigma'$ are 
non-singular.  Let $X'$ be the toric variety corresponding to $\Sigma'$.  
Then $X'$ is smooth, and there is a proper birational map from $X'$ to $X$. 
The Todd class of $X$ coincides with the push-forward of the Todd class of 
$X'$.
Supposing that we have a complement map $\Psi$
for $X'$, we can thus define a cycle 
Todd class for $X$, written $\ctd^{X'}(X)$, as the pushforward of 
$t_\Psi(X')$.

A priori, this cycle Todd class $\ctd^{X'}(X)$ depends
 on the choice of subdivision 
$\Sigma'$ of $\Sigma$.  
We will show that under reasonable assumptions, $\ctd^{X'}(X)$ does not
depend on $\Sigma'$.  
The chief result of this paper is the following theorem and its corollary:

\begin{theorem} Let $X$ be any toric variety, with fan $\Sigma$.  
Let $\Psi$ be a rigid complement map, and $X'$ the 
toric variety induced by any non-singular subdivision $\Sigma'$ of 
$\Sigma$ such that $\Psi$ is a complement map for $\Sigma'$.  Then the cycle 
Todd class $\ctd^{X'}(X)$ does not depend on the choice of $X'$.  \end{theorem}
Before giving the proof, we note that Theorem 1 allows us to
construct from any rigid complement map,  a cycle Todd class
with rational coefficients that satisfies
Danilov's criterion.

\begin{corollary} Fix a rigid complement map $\Psi$.
\begin{enumerate}
\item[(i)] Let
$\sigma$ be a cone in the domain of $\Psi$.  Then, for any fan $\Sigma$
containing $\sigma$, the coefficient of $[V(\sigma)]$ in the cycle Todd 
class of the
toric variety defined by $\Sigma$ depends only on the cone $\sigma$
(and is otherwise independent of $\Sigma$).  We 
call this coefficient the Todd measure of $\sigma$, and we
denote it $\mu(\sigma)$.  

\item[(ii)] {\it (Danilov's question)} If $\Psi$ is a rigid complement map
 defined on all cones in $N$, then the
 function $\mu$ obtained in this way
satisfies Danilov's criterion: for any fan $\Sigma$ in $N$ with associated
toric variety $X$, we have
$$
\operatorname{Td}\, X = \sum \mu(\sigma) [V(\sigma)].
$$

\item[(iii)]  {\it (Additivity)} If $\sigma_1,\dots,\sigma_t$ are cones of 
the same dimension which overlap only on boundaries, such that 
their union is also a (strictly convex) cone $\tau$, and
all the $\sigma_i$ and $\tau$ are in the domain of $\Psi$, then
$$ \mu(\tau)=\sum_{i=1}^t \mu(\sigma_i).$$ 
(Note that by the rigidity of $\Psi$, if any of the $\sigma_i$ or $\tau$
is in the domain of $\Psi$, so are all the others.)  

\item[(iv)]
Let $P$ be a lattice polytope in the dual lattice to $N$,
$M=\Hom(N,\Bbb{Z})$.   For
$F$ a face of $P$, let $C(F,P)$ be the outer normal cone to $F$ in $N$.  
Let
$\lvol(F)$ denote the volume of $F$, normalized with respect to the lattice
in the affine span of $F$.  

Suppose that the outer normal fan of $P$ has a non-singular subdivision all
of whose cones are in the domain of $\Psi$.    
Then we have the following formula for the
number of lattice points of $P$:
$$\#(P)=\sum_{F} \mu(C(F,P))\lvol(F),$$
where the sum is taken over all faces $F$ of $P$.
\end{enumerate}
\end{corollary}

\begin{proof} \begin{enumerate} \item[(i)]
We fix one way to refine $\sigma$ to a smooth fan, and 
then, using that refinement, we will always obtain the same coefficient 
of $[V(\sigma)]$ independent of $\Sigma$. By Theorem 1, this coefficient
is also independent of the choice of refinement.  

\item [(ii)] This follows immediately from the definition of $\mu$.

\item[(iii)]  Let $\Sigma$ be a fan containing $\tau$, and $X$ the associated
toric variety.  Let $\Sigma'$ be the fan
obtained by replacing $\tau$ by $\sigma_1,\dots,\sigma_t$, and $X'$ the
associated toric variety.  Let $Y\rightarrow X'$ be a resolution of 
singularities of $X'$.  The desired result follows from the fact that
the coefficient of $[V(\tau)]$ in $\ctd^Y(X)$ is the sum of the 
coefficients of $[V(\sigma_i)]$ in $\ctd^Y(X')$.  

\item [(iv)] 
This follows from (ii) by an application of Hirzebruch-Riemann-Roch,
see [Dan] or [Ful].  
\end{enumerate} \end{proof}

To prove Theorem 1, we shall connect the cycle Todd classes 
$\ctd^{X'}(X)$ to the 
Todd class formula of \cite{Pom2}.

\section{A second approach to cycle Todd classes}

The Todd class formula from \cite{Pom2} expresses the Todd class of a 
complete, simplicial toric variety in a canonical way 
as a product of Cartier divisors: for complete simplicial $X$, a polynomial 
$p_X$ is
defined such that $p_X(D_1,\dots,D_r)$, evaluated in the Chow ring with
rational coefficients, gives the Todd class.  Fix a complement map
$\Psi$ for $X$. 
We define $\ctd^s(X)$ to be
the result of evaluating $p_X(D_1,\dots,D_r)$ in $Z_*(X)_\Q$.  
It follows that $\ctd^s(X)$ is a cycle Todd class of $X$.  

We prove the following theorem:

\begin{theorem} For $X$ a complete, simplicial toric variety, 
$X'\rightarrow X$ a toric resolution of singularities of $X$, and 
$\Psi$ a 
rigid complement map for $X'$, 
then $\ctd^s(X)=\ctd^{X'}(X)$. \end{theorem}

\begin{proof}
We begin by recalling the notion of a stellar subdivision of a fan.  
Let $\Sigma$ be a fan, and $\rho_0$ a ray which is not a ray of $\Sigma$.
Let $\sigma$ be the minimal cone of $\Sigma$ containing $\rho_0$.  
The stellar subdivision of $\Sigma$ at $\rho_0$ consists of replacing 
every cone $\tau$ which contains $\rho_0$ by the cones $\rho_0+\delta$
where $\delta$ is a face of $\tau$ not containing $\rho_0$.  

The main element in the proof of Theorem 2 is the following special
case of the theorem:
  
\begin{lemma} Let $\pi: Y\rightarrow Z$ be a proper birational map of
simplicial toric varieties corresponding to a stellar subdivision.  Then

$$\ctd^{Y}(Z)=\ctd^s(Z).$$\end{lemma}

\begin{proof}
Let $\rho_0$ be the ray that appears in the fan of $Y$ and not $Z$.  
Let $D_0$ be the 
$\Q$-Cartier divisor on $Y$ corresponding to $V(\rho_0)$.
Let
$\sigma$ be the minimal cone of the fan of $Z$ containing $\rho_0$.  
Let $E_1,\dots,E_d$ be the $\Q$-Cartier divisors on $Z$ corresponding to rays
in $\sigma$.  Let $E_{d+1},\dots,E_r$ be the $\Q$-Cartier divisors on $Z$ 
corresponding to rays not in $\sigma$. For $1\leq i \leq r$, 
let $D_i$ be the $\Q$-Cartier divisor
on $Y$ corresponding to the same ray as $E_i$.

In \cite{Pom2}, a map $\tilde \pi$ from $\Q[y_0,\dots,y_r]$ to 
$\Q[z_1,\dots,z_r]$ 
is defined as follows:
\begin{enumerate}

\item[A)] $\tilde \pi(y_if)=z_i\tilde \pi(f)$ for $i>d$

\item[B)] $\tilde \pi(y_if)=\frac{m_i}{m_j}\tilde \pi (y_jf) +( z_i -
\frac{m_i}{m_j}z_j)\tilde \pi(f)$ for $1\leq i,j \leq d$

\item[C)] $\tilde \pi(y_0f)=-\frac{m_0}{m_j}\tilde\pi(y_jf) + \frac{m_0}{m_j}
z_j\tilde\pi (f)$ for $1\leq j\leq d$

\item[D)] $\tilde \pi(y_1\dots y_df)=0$

\item[E)] $\tilde \pi (g(y_1,\dots,y_r))=g(z_1,\dots,z_r)$ for $g$ 
any polynomial of degree at most 
$d-1$.  
\end{enumerate}
Here 
 $m_0$ is the multiplicity of 
$\sigma$ (i.e.~the index of the lattice generated by the first lattice
points along the extreme rays of $\sigma$ in the intersection of 
$N$ with the span of $\sigma$), and $m_i$ is 
the multiplicity of the cone generated by the rays corresponding to 
$D_0,\dots,\hat D_i, \dots, D_r$.  

It is relatively easy to show that there is at most one map 
$\tilde \pi$ with these properties.  We will review the argument
below, as we shall have need of it in a slightly different context.  
It is much harder to show that there is a map $\tilde\pi$
with all these properties.  
The reader is referred to \cite{Pom2} for the proof of this. 

By construction, $p_Z=\tilde \pi(p_Y)$. Therefore, we may apply
the following lemma to show 
that $p_Z(E_1,\dots,E_r)$ evaluated in $Z_*(Z)_\Q $ coincides with
the pushforward of $p_Y(D_0,\dots,D_r)$ evaluated in $Z_*(Y)_\Q$, which
will prove Lemma 2.  

\begin{lemma} Let $f$ be a polynomial in $y_0,\dots,y_r$.  Then we
have the following equality in $Z_*(Z)_\Q$:
$$\tilde\pi(f)(E_1\dots,E_r)= \pi_*(f(D_0,\dots,D_r)).$$  
\end{lemma}

\begin{proof}  
First, we check the following properties, 
where $W$ is any cycle in $Z_*(Y)$ and the $m_i$ are as above:

\begin{enumerate}
\item[A$'$)] $\pi_*(D_i W)=E_i(\pi_*(W))$ for $i>d$

\item
[B$'$)] $\pi_*(D_iW)= \frac{m_i}{m_j}\pi_* (D_jW) +( E_i -
\frac{m_i}{m_j}E_j)\pi_*(W)$ for $1\leq i,j \leq d$

\item
[C$'$)] $\pi_*(D_0W)=-\frac{m_0}{m_j}\pi_*(D_jW) + \frac{m_0}{m_j}
E_j\pi_* (W)$ for $1\leq j\leq d$

\item
[D$'$)] $\pi_*(D_1\dots D_dW)=0$

\item[E$'$)] $\pi_* (g(D_1,\dots,D_r))=g(E_1,\dots,E_r)$ for $g$ 
any polynomial of degree at most 
$d-1$.  

\end{enumerate}

For $i>d$, $\pi^*(E_i) = D_i$.  Thus, by an application of Lemma 1,
$\pi_*(D_iW)=E_i\pi_*(W)$.  This establishes A$'$.  B$'$ and C$'$ follow
similarly, using the facts that $\pi^*(m_jE_i - m_iE_j)=m_jD_i-m_iD_j$
and $\pi^*(m_0E_j)=m_0D_j+m_jD_0$.  D$'$ follows because $D_1\dots D_d=
0 \in Z_*(Y)_\Q$, since there are no cones containing all the corresponding
rays.  To establish E$'$, we observe that 
$g(D_1,\dots,D_r)$ evaluated in $Z_*(Y)_\Q$ will consist of a sum of 
$V(\tau)$ with $\tau$ of codimension $d-1$ or less.  
For those $\tau$ not containing 
the ray corresponding to $D_0$, $V(\tau)$ will occur with the same coefficient
in $g(D_1,\dots,D_r)$ and in $g(E_1,\dots,E_r)$.  On the other hand, if
$\tau$ does contain this ray, then $\pi_*([V(\tau)])=0$, because the image
of $V(\tau)$ is $V(\sigma)$, which has lower dimension.  

Next, we use essentially the argument of \cite{Pom2} which establishes
the uniqueness of $\tilde \pi$.  Consider a monomial $Q=D_0^{a_0}
\dots D_r^{a_r}$.
We wish to show that $\pi_*(D_0^{a_0}\dots D_r^{a_r})$ is determined
by properties A$'$) through E$'$).  The proof is inductive.  E$'$) 
establishes the 
result for terms of degree at most $d-1$.  Now assume the result for terms
of lower degree than $Q$.      
If any $a_i$ are 
greater than zero for $i>d$, apply A$'$), and the desired result follows by
induction.  
If there is some $1\leq j \leq d$ such that $a_j$ is zero, apply 
either B$'$) or C$'$) to 
replace $Q$ by terms which can be dealt with by induction,
plus  a term which has fewer $i$, $1 \leq i \leq d$, such that $a_i=0$.  
Repeat this process until the problematic term has $a_i \ne 0$ for 
$1\leq i \leq d$, at which point it can be disposed of by D$'$).  

Now observe that, by properties A) through E), the map taking a polynomial
$f\in \Q[y_0,\dots,y_r]$ to  
$\tilde \pi(f)(E_1,\dots,E_r)\in Z_*(Z)_\Q$ has properties A$'$) through 
E$'$), and 
thus
it coincides with the map taking $f$ to $\pi_*(f(D_0,\dots,D_r))$, 
as desired.  

This establishes Lemma 3. 
As already described, Lemma 2 follows.  \end{proof} \end{proof}

Now, we turn to the proof of Theorem 2.  
Let $W=W_p\rightarrow W_{p-1}\rightarrow \dots\rightarrow W_0=X$ be a 
resolution of singularities for $X$ such that the map from $W_i$ to 
$W_{i-1}$ is induced by a stellar subdivision of fans.  (The existence
of a resolution of singularities of this type for any $X$ is shown in
[Ful].)  Let $w_i$ denote the map from $W_i$ to $X$.  

Now, we prove by induction on $i$ that $(w_i)_*(\ctd^s(W_i))=\ctd^s(X)$.
The $i=0$ case is trivial.  The induction step consists in observing that
$w_{i+1}$ factors through $W_i$, and then applying Lemma 2.  
This establishes that $\ctd^W(X)=\ctd^s(X)$.

We must now prove that $\ctd^W(X)=\ctd^{X'}(X)$ for $X' \rightarrow
X$ any resolution of singularities.  
We recall the Weak Oda Theorem:

\begin{kludge} 
Given $\Pi$ and $\Sigma$ two non-singular fans
with the same support, there is a sequence of non-singular 
fans $\Pi=\Delta_0,
\Delta_1,\dots,\Delta_p=\Sigma$ such that either $\Delta_{i+1}$ is
obtained from $\Delta_{i}$ by a stellar subdivision or vice versa.  

Further, if both $\Pi$ and $\Sigma$ refine some fan $\Omega$, then all the
$\Delta_i$ may be chosen so as to refine $\Omega$.  
\end{kludge}

Applying the Weak Oda Theorem to $W$ and $X'$ (both of whose fans refine
the fan for $X$), we obtain a 
sequence of toric varieties $X'=Z_p,\dots,Z_0=W$, all of whose fans refine
the fan of $X$, such that either $Z_{i+1}$ is obtained from $Z_i$ by a 
stellar subdivision or vice versa.  Lemma 2 implies that
$\ctd^{W_i}(X)=\ctd^{W_{i+1}}(X)$, and it follows that $\ctd^{X'}(X)=
\ctd^{W}(X)$, as desired.  
\end{proof}

We now prove Theorem 1.

\begin{proof}  
If $X$ is complete and simplicial, then by Theorem 2, $\ctd^{X'}(X)=
\ctd^s(X)$, so $\ctd^{X'}(X)$ does not depend on $X'$.  However, we must establish the 
result without the assumption that $X$ is complete and simplicial.  

Suppose first that $X$ is complete but not simplicial.  
Let $Y\rightarrow X$ and $Y' \rightarrow X$ be two 
resolutions of singularities.  We can always construct a resolution of 
singularities $Z\rightarrow X$ which factors through both $Y$ and $Y'$,
by taking a fan which is a subdivision of both $Y$ and $Y'$, and then
subdividing it further so as to make it non-singular.  
Thus, it suffices to show that if $g:Z \rightarrow Y$ and $f:Y \rightarrow
X$, such that $f$ and $f\circ g$ are resolutions of singularities for $X$,
then $\ctd^Y(X)=\ctd^Z(X)$.  Now 
\begin{eqnarray*}
\ctd^Z(X)&=&(f\circ g)_* (\ctd(Z))\\
&=&f_*(g_*(\ctd(Z)))\\
&=&f_*(\ctd(Y)) \\
&=&\ctd^Y(X)
\end{eqnarray*}
where we have used the fact that $g_*(\ctd(Z))=\ctd(Y)$, by Theorem 2.

Now suppose that $X$ is not complete.  Let $\tilde X$ be a toric variety
obtained by adding cones to the fan $\Sigma$ for $X$, in order to make it
complete.  As before, let $Y\rightarrow X$ and $Y' \rightarrow X$ be 
resolutions of singularities, and let $\tilde Y$ and $\tilde Y'$ be 
resolutions of singularities for $\tilde X$ whose fans intersected with
the support of $\Sigma$ 
agree with those of $Y$ and
$Y'$ respectively. For
any $\sigma\in\Sigma$, the coefficient of $[V(\sigma)]$ in $t^Y_\Psi(X)$
is the same as the coefficient of $t^{\tilde Y}_\Psi(\tilde X)$, and similarly
replacing $Y$ by $Y'$ and $\tilde Y$ by $\tilde Y'$.  Thus, the general
result follows by applying the result in the complete case to $\tilde X$.  
This proves Theorem 1.    \end{proof}

\section{Cycle Todd classes induced by a choice of flag}

In this section we investigate further $\ctd(X)$ for $X$ smooth and 
$\Psi$ induced by a complete flag. Our tool will be a more general
result about computation in $Z_*(X)_\Q$ when $\Psi$ is induced from a 
complete flag.

\begin{theorem}Let $X$ be a simplicial toric variety, and
let $\Psi$ be induced from a complete flag $F_\bullet$. Let $k\leq n$,
and let  $a_1,\dots, a_k$ be non-negative integers that sum to $k$.
Let $\sigma$ be the cone spanned by $\rho_1,\dots,\rho_k$. 
Let $v_i$ be the first lattice point
along $\rho_i$.  By the genericity of $F_\bullet$, there is a unique
(up to scalar multiple) collection of non-zero $t_i$ such that
$$\sum_{i=1}^k t_iv_i\in F_{n-k+1}.$$
Then the coefficient of $[V(\sigma)]$ in 
$\prod_{i=1}^kD_i^{a_i}$ is 
$$\frac{1}{\mult\sigma} \prod_{i=1}^k t_i^{a_i-1}.$$
\end{theorem}

\begin{proof}
Most of our effort will be directed to showing that the coefficient of 
$[V(\sigma)]$ in $\prod_{i=1}^kD_i^{a_i}$ coincides with the 
coefficient of $[V(\sigma)]$ in 
$\prod_{i=1}^k t_i^{a_i-1}D_i$.  
The proof is by induction on the number of $a_i$ which equal zero.  The
base case, when none of the $a_i=0$, and thus all the $a_i=1$, 
is tautological.

So suppose there is some $a_l=0$.  Pick $j$ such that $a_j>1$.  
Observe that there is an element $m\in F_{n-k+1}^\perp$
such that $\lp m,v_j\rp=t_l$,
$\lp m,v_l\rp=-t_j$, and $\lp m,v_i\rp=0$ for all $i\ne j,l$.   
By Proposition 3, 
for any $\tau\prec\sigma\in \Sigma$ of dimension $k-1$,
$$
\div(\chi^m)\cdot [V(\tau)]=0
$$
and thus the coefficients of $[V(\sigma)]$ in 
$t_l D_j\cdot [V(\tau)]$ and $t_j D_l\cdot [V(\tau)]$ coincide.
  
Consequently, the coefficient of $[V(\sigma)]$ in $D_1^{a_1}\cdots D_k^{a_k}$
is the same as the coefficient of $[V(\sigma)]$ in 
$$
\frac {t_j}{t_l} D_1^{a_1}\cdots D_j^{a_j-1}\cdots D_l^{a_l+1}\cdots D_k^{a_k}.
$$
This expression has one fewer $a_i$ which is zero, and our claim follows
by induction.  

The statement of the corollary now follows from the fact that 
$$D_1\cdot\cdots\cdot D_k=\frac{1}{\mult \sigma} [V(\sigma)].$$  
\end{proof}

The following corollary states that for a complement map induced
by a complete flag, the Todd class $t_{\Psi}(X)$ we have constructed
coincides with the Todd class of \cite[Theorem 11]{Mor}. For
each $k$-dimensional cone $\sigma$, Morelli
defines $\mu_k^{td_k}(\sigma)$, a rational function on the 
Grassmannian $Gr_{n-k+1}(N)$ of $n-k+1$-planes in $N$.  In our 
$t_{\Psi}(X)$, the coefficient $\mu(\sigma)$ of the torus-invariant cycle
$[V(\sigma)]$ depends only 
on $F_{n-k+1}$, 
where $k$ is the dimension of $\sigma$, and not on the rest of the flag.
In this way we can view $\mu$ as 
a rational function on $Gr_{n-k+1}(N)$.  The corollary
below assures us that this function is exactly  $\mu_k^{td_k}(\sigma)$.

\begin{corollary} Let  $\mu_k^{td_k}$ denote the function
of \cite[Theorem 11]{Mor}.  Let $F$ be a complete flag, let
$\Psi$ be the induced complement map, and let $\mu$ be the
associated Todd measure (as in Corollary 1(i)).  Then for any
cone of dimension $k$,
$$
\mu(\sigma) = \mu_k^{td_k}(\sigma)(F_{n-k+1}).
$$
\end{corollary}

\begin{proof}
For nonsingular $\sigma$, the equation follows by comparing Theorem 3
to the definition of $\mu_k^{td_k}$ \cite[p. 218]{Mor}.

Suppose that $\sigma$ is singular.  Let $\Sigma$ be a fan consisting
of $\sigma$ and its faces, and let $\Sigma'$ be a desingularization of 
$\Sigma$.  Let $\tau_1,\dots,\tau_m$ be the $k$-dimensional cones of 
$\Sigma'$.  By the additivity of $\mu_k^{td_k}$ \cite[p. 184]{Mor},
$\mu_k^{td_k}(\sigma)=\sum_{i=1}^m \mu_k^{td_k}(\tau_i)$.  By the 
additivity of $\mu$ (Corollary 1(iii)), the same formula holds with
$\mu$ replacing $\mu_k^{td_k}$ everywhere, which reduces us to the
non-singular case.  
\end{proof}

\section{Canonical cycle Todd classes}

A drawback of the cycle Todd classes presented so far is that they depend
on a choice of complement map.  Further, no single choice of fan defines
a complement map simultaneously for all cones in $N$.  Both these problems
can be avoided, as is done in \cite{Mor}, by enlarging the coefficient field
of the cycles.  In this way we obtain a canonical Todd class based on the
inner product construction, and one based on the complete flag construction.

Let $\Cal{F}$ be the set of all complete flags in $N_\Q$, which has the
structure of an algebraic variety. For $F$ a complete flag in $N_\Q$,
generic with respect to the fan for $X$, let $t_F(X)$ be the
cycle Todd class for $X$ given by the complement map induced by $F$.
Now, define $t_\Cal{F}(X):\Cal{F}\dottedrightarrow Z_*(X)_\Q$ by 
$t_\Cal{F}(X)(F)=t_F(X)$ when $t_F(X)$ is defined.  We observe that, by the 
results of the previous section,
$t_\Cal{F}(X)$
is a rational function on $\Cal{F}$, and thus if we write $\Bbb{K}$ for the
rational function field of $\Cal{F}$, we can view $t_\Cal{F}(X)$ as lying in
$Z_*(X)_\Bbb{K}$, and it is clear that it is a cycle Todd class for 
$X$.  

One can do essentially the same thing to produce a canonical Todd class
based on the inner product construction.  
We construct Todd classes in the rational
function field of $\Cal{S}$, the set of all symmetric $n\times n$ matrices. 
Fixing a basis for $N$, we can
identify inner products on $N_\Q$ with symmetric positive definite
$n\times n$ matrices.  Using the approach of the previous paragraph,
we define cycle Todd classes whose coefficients are functions from 
symmetric positive definite $n\times n$ matrices to $\Q$; one then
checks that these functions extend to rational functions on $\Cal{S}$.

\section*{Acknowledgements}

The authors thank William Fulton and Burt Totaro for much useful advice,
and thank the referee for many helpful suggestions, including an improvement
to the statement of Theorem 3.

\end{document}